\pdfoutput=1
\RequirePackage{ifpdf}
\ifpdf % We are running pdfTeX in pdf mode
\documentclass[pdftex]{sigma}
\else
\documentclass{sigma}
\fi

\numberwithin{equation}{section}

\newtheorem{Theorem}{Theorem}[section]
\newtheorem{Corollary}[Theorem]{Corollary}
\newtheorem{Lemma}[Theorem]{Lemma}
\newtheorem{Proposition}[Theorem]{Proposition}
 { \theoremstyle{definition}
\newtheorem{Definition}[Theorem]{Definition}
\newtheorem{Example}[Theorem]{Example}
\newtheorem{Remark}[Theorem]{Remark} }

\begin{document}

\allowdisplaybreaks

\renewcommand{\thefootnote}{$\star$}

\renewcommand{\PaperNumber}{068}

\FirstPageHeading

\ShortArticleName{Dynamics on Wild Character Varieties}

\ArticleName{Dynamics on Wild Character Varieties\footnote{This paper is a~contribution to the Special Issue on
Algebraic Methods in Dynamical Systems.
The full collection is available at
\href{http://www.emis.de/journals/SIGMA/AMDS2014.html}{http://www.emis.de/journals/SIGMA/AMDS2014.html}}}

\Author{Emmanuel PAUL~$^\dag$ and Jean-Pierre RAMIS~$^\ddag$}

\AuthorNameForHeading{E.~Paul and J.-P.~Ramis}

\Address{$^\dag$~Institut de Math\'ematiques de Toulouse, CNRS UMR 5219, \'Equipe \'Emile Picard,\\
\hphantom{$^\dag$}~Universit\'e Paul Sabatier (Toulouse 3), 118 route de Narbonne,\\
\hphantom{$^\dag$}~31062 Toulouse CEDEX 9, France}
\EmailD{\href{mailto:emmanuel.paul@math.univ-toulouse.fr}{emmanuel.paul@math.univ-toulouse.fr}}

\Address{$^\ddag$~Institut de France (Acad\'emie des Sciences) and Institut de Math\'ematiques de Toulouse,\\
\hphantom{$^\ddag$}~CNRS UMR 5219, \'Equipe \'Emile Picard,  Universit\'e Paul Sabatier (Toulouse 3),\\
\hphantom{$^\ddag$}~118 route de Narbonne, 31062 Toulouse CEDEX 9, France}
\EmailD{\href{mailto:ramis.jean-pierre@wanadoo.fr}{ramis.jean-pierre@wanadoo.fr}}

\ArticleDates{Received March 26, 2014, in f\/inal form August 05, 2015; Published online August 13, 2015}

\Abstract{In the present paper, we will f\/irst present brief\/ly a general research program about the study of the ``natural dynamics'' on character varieties and wild character va\-rieties. Afterwards, we will illustrate this program in the context of the Painlev\'e dif\/ferential equations~$P_{\rm VI}$ and~$P_{\rm V}$.}


\Classification{34M40; 34M55}

\begin{flushright}
\emph{To Juan J. Morales-Ruiz, for his $60^{th}$ birthday.}
\end{flushright}

\renewcommand{\thefootnote}{\arabic{footnote}}
\setcounter{footnote}{0}

\section{A sketch of a program}

We begin with the sketch of a \emph{work in progress} of the authors in collaboration with Julio Rebelo, based on (or related to) some results due to several people, mainly:
Ph.~Boalch \cite{Bo1, Bo2, Bo3, Bo4, Bo5}, S.~Cantat, F.~Loray~\cite{CanLo}, B.~Dubrovin~\cite{Du,DuMaz}, M.A.~Inaba, K.~Iwasaki~\cite{IWA}, M.~Jimbo~\cite{JMU}, B.~Malgrange~\cite{Ma1, Ma2, Ma4, Ma5, Ma3}, M.~Mazzoco, T.~Miwa, M.~van der Put, M.-H.~Saito~\cite{VdPSai},  K.~Ueno~\cite{Ue}, E.~Witten~\cite{Wi}, \dots,
and the Kyoto school around T.~Kawai and Y.~Takei~\cite{KaKNT, KaT}.

In the present state it is mainly a  \emph{PROGRAM}.

We would like to understand:
\begin{enumerate}\itemsep=0pt
\item
The \emph{dynamics} and the
 \emph{wild dynamics}\footnote{That is, roughly speaking, the ordinary dynamics coming from the nonlinear monodromy ``plus'' the dynamics coming from ``nonlinear Stokes phenomena''.} of equations of  \emph{isomonodromic deformations} and of \emph{wild isomonodromic deformations} using the (\emph{generalized\,})
 \emph{Riemann--Hilbert correspondances} and the corresponding (wild) dynamics on the  \emph{$($wild$)$ character varieties}. The notion of wild character variety was introduced by Boalch. The braid group action on character varieties for the Painlev\'e equations has been f\/irst def\/ined by Dubrovin and Mazzocco for special parameters in~\cite{DuMaz}, and by Iwasaki~\cite{IWA} in the general case.
\item
The \emph{confluence phenomena} for the equations of (wild) isomonodromic deformations and the corresponding conf\/luence phenomena for the (wild) dynamics.
\end{enumerate}

Our (long term!) aim is to built a \emph{general theory}, testing it at each step on the case of the \emph{Painlev\'e equations} (which is already far to be trivial).

Our initial motivation was to compute the \emph{Malgrange groupoids}
of the six Painlev\'e equations. Our conjecture is that it is  \emph{the biggest possible} (that is the groupoid of transformations conserving the area)
in the generic cases (it is known for $P_{\rm I}$: see {\small Casale} in \cite{CAS}, for $P_{\rm VI}$: Cantat--Loray in~\cite{CanLo}, and for special parameters of $P_{\rm II}$ and $P_{\rm III}$: Casale and Weil~\cite{CaWe}).
Using a result of Casale~\cite{CAS}, it is possible to reprove, in ``Painlev\'e style'', the
\emph{irreducibility}
of Painlev\'e equations (initially proved by the japanese school: Nishioka, Umemura, Okamoto, Noumi, \dots). Our approach is to try to def\/ine in each case, using the Riemann--Hilbert map, a \emph{wild dynamics} on the Okamoto variety of initial conditions and to prove that this dynamics is ``chaotic'', forcing Malgrange groupoid to be big (up to the conjecture that the wild dynamics is ``into'' the Malgrange groupoid).

The classical character varieties are moduli spaces of monodromy data of regular-singular connections, that is spaces of
representations of the fundamental group of a punctured (or not) Riemann surface. Atiyah-Bott and Goldman prove that they admit holomorphic Poisson structure. This fact has been extended to wild character varieties by Boalch (see~\cite{Bo0,Bo4}).

The wild character varieties generalize the classical (or \emph{tame}) character varieties.
They are moduli spaces of generalized monodromy data of meromorphic connections. In the irregular case it is necessary to add ``Stokes data" to the classical monodromy. Then the wild character varieties are spaces of representations of a wild fundamental
\emph{groupoid}.

In the global irregular case it is \emph{necessary} to use groupoids. They are explicitely used in \cite{Bo0}, and implicitely used in \cite{JMU, Wi}.
In the local irregular case it is suf\/f\/icient to use a group, the Ramis wild fundamental group~\cite{MR2, VdPSin}.

Therefore in order to understand the conf\/luence process of a classical representation of the fundamental group towards a representation of the wild groupoid, it is better to
replace the classical fundamental group by a groupoid. This is a posteriori clear in the computations of~\cite{Ra1} in the hypergeometric case. We plan to return to such problems in future papers.

We will show below, with the example of~$P_{\rm VI}$, that even in the classical case
it is better to use fundamental groupoids than fundamental groups to study character varieties and their \emph{natural dynamics}. This is in a line strongly suggested by Alexander Grothendieck.

\emph{Ceci est li\'e notamment au
fait que les gens s'obstinent encore, en calculant avec des groupes fondamentaux,
\`a fixer un seul point base, plut\^ot que d'en choisir astucieusement tout
un paquet qui soit invariant par les sym\'etries de la situation, lesquelles sont
donc perdues en route. Dans certaines situations (comme des th\'eor\`emes
de descente \`a la Van Kampen pour les groupes fondamentaux) il est bien plus
\'el\'egant, voire indispensable pour y comprendre quelque chose, de travailler
avec des groupo\"\i des fondamentaux par rapport \`a un paquet de points base
convenable, et il en est certainement ainsi pour la tour de Teichm\"uller} (cf.~\cite{Gro}).

\emph{\dots people are accustomed to work with fundamental groups and generators and relations for these and stick to it, even in contexts when this is wholly inadequate, namely when you get a clear description by generators and relations only when working simultaneously with a whole bunch of base-points chosen with care~-- or equivalently working in the algebraic context of groupoids, rather than groups. Choosing paths for connecting the base points natural to the situation to one among them, and reducing the groupoid to a single  group, will then hopelessly destroy the structure and inner symmetries of the situation, and result in a mess of generators and relations no one dares to write down, because everyone feels they won't be of any use whatever, and just confuse the picture rather than clarify it. I have known such perplexity myself a long time ago, namely in Van Kampen type situations, whose only understandable formulation is in terms of $($amalgamated sums of$)$ groupoids} (Alexandre Grothendieck, quoted by Ronald Brown\footnote{\url{http://pages.bangor.ac.uk/~mas010/pstacks.htm}.}).

We recall that a groupoid is a small category in which every morphism is an
isomorphism (for basic def\/initions and details cf.~\cite{Br}).

\begin{Example} Let $Y$ be a (topological) manifold. The fundamental groupoid $\pi_1(Y)$ of $Y$  is the groupoid whose objects are the elements~$y$ of $Y$, and whose morphisms are the paths between elements of~$Y$ up to homotopy.
\end{Example}

We have the following generalization.

\begin{Definition}
Let $Y$ be a (topological) manifold and $S\subset Y$. The fundamental groupoid $\pi_1(Y,S)$ of the pair $(Y,S)$  is the groupoid whose objects are the elements $s$ of $S$ in $Y$, and whose morphisms are the paths between elements of $S$ up to homotopy.
\end{Definition}

We have $\pi_1(Y)=\pi_1(Y,Y)$ and when $S$ is reduced to a point $a$, $\pi_1(Y,\{a\})$ is the classical fundamental group of $Y$ based at $a$: $\pi_1(Y,\{a\})=\pi_1(Y,a)$.

\section[Character varieties in the Painlev\'e context: the regular singular case revisited]{Character varieties in the Painlev\'e context: \\ the regular singular case revisited}

We focus now on the ``Painlev\'e context'':
let $\mathcal{E}$ be the set of linear rank~2 connections on the trivial bundle over~$\mathbb{P}_1(\mathbb{C})$, with coef\/f\/icients in ${\mathfrak{sl}}_2(\mathbb{C})$ and such that its singular locus contains at most 4 singular points. In this section, we f\/irst consider the ``classical'' case with~$4$ \emph{regular} singular points, in order to be more familiar with the groupoid point of view which is essential to deal with the irregular cases. Furthermore it turns out that this point of view is yet useful to obtain the dynamics in the regular singular case.
In this section, $\Delta$~is a linear dif\/ferential system which represents the connection $\nabla$.

\subsection{The fundamental groupoid}
We consider the ``extended'' singular locus $S$ of~$\Delta$: $S$ is the set of pairs $s=(p,d)$ where $p$ is a~singular point of~$\Delta$ and $d$ is a~ray based in $p$. Therefore $s$ is also a point on the divisor~$D_p$ of the real blowing up~$E_p$ at~$p$.
Let $X$ be the manifold obtained by the real blowing up of each singular point. We denote by $\gamma_{s,s}$ a loop from~$s$ to~$s$ in~$X$ with positive orientation, homotopic to the exceptional divisor:

\begin{figure}[h!]\centering
 \includegraphics[scale=0.6]{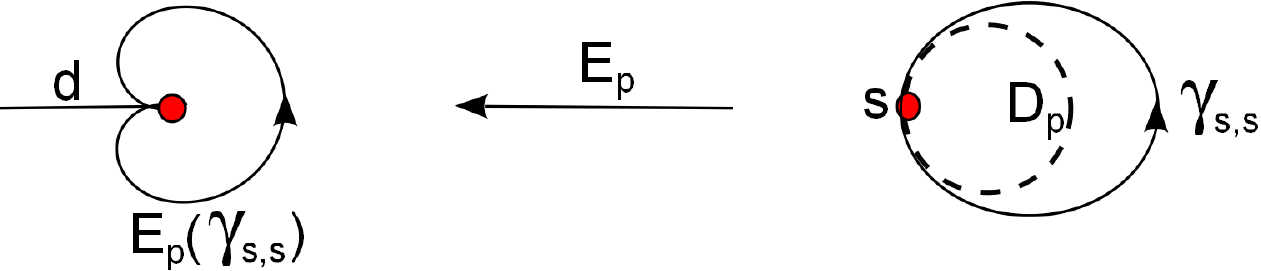}
\caption{The real blowing up at $p$.}\label{figure1}
\end{figure}

\begin{Definition}
The fundamental groupoid $\pi_1(X,S)$  is the groupoid whose objects are the elements $s$ of $S$ in $X$, and whose morphisms are the paths between elements of $S$ up to homotopy.
The subgroupoid $\pi_1^{\rm loc}(X,S)$ is the groupoid with same objects, whose morphisms are generated only by the loops~$\gamma_{i,i}$ homotopic to each exceptional divisor at~$p_i$.
\end{Definition}

We denote:
\begin{enumerate}\itemsep=0pt
\item[--] $\operatorname{Aut}(\pi_1(X,S))$ the group of the automorphisms of the groupoid, and $\operatorname{Aut}_0(\pi_1(X,S))$ the subgroup of the ``pure'' automorphisms, which f\/ix each object.

\item[--] $\operatorname{Inn}_0(\pi_1(X,S))$ the normal subgroup of $\operatorname{Aut}_0(\pi_1(X,S))$ of the inner automorphisms. An inner automorphism is def\/ined by a collection of loops $\alpha_{i}$ at each object, by setting:
    \begin{gather*}
    h_{\{\alpha_i\}}\colon \ \gamma_{i,j}\mapsto  \alpha_i\gamma_{i,j}\alpha_j^{-1} .
    \end{gather*}
\item[--] $\operatorname{Out}_0(\pi_1(X,S)):=\operatorname{Aut}_0(\pi_1(X,S))/\operatorname{Inn}_0(\pi_1(X,S))$.

\item[--] $\operatorname{Out}_0^*(\pi_1(X,S))$:= the subgroup of $\operatorname{Out}_0(\pi_1(X,S))$ whose elements f\/ix each local morphism $\gamma_{i,i}$ up to conjugation.
\end{enumerate}

We obtain the following presentation of $\pi_1(X,S)$:

\begin{figure}[h!]
\includegraphics[scale=0.5]{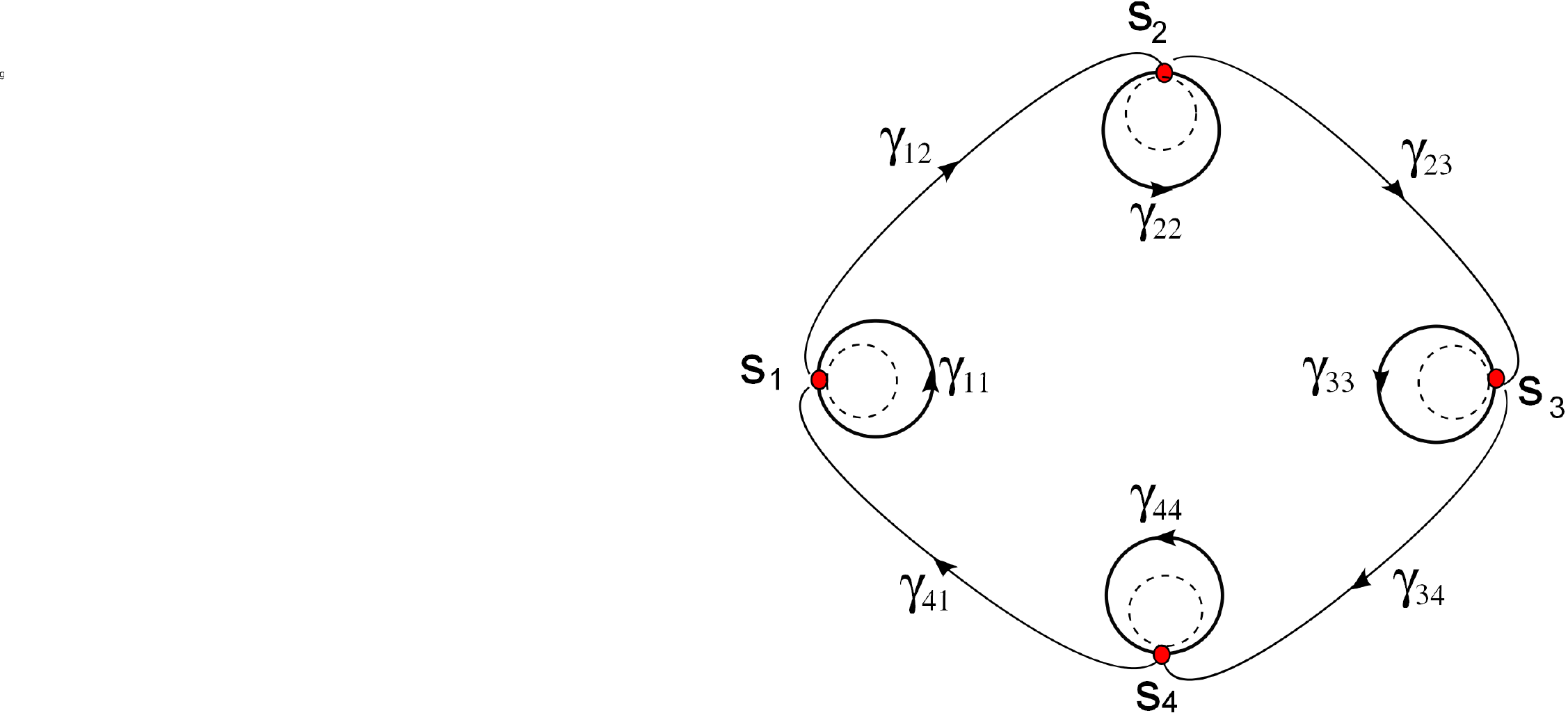}
\caption{The groupoid $\pi_1(X,S)$ (classical case).}\label{figure1bis}
\end{figure}

The exceptional divisors are circles in dotted lines. The morphisms are generated by 8 paths~$\gamma_{i,i}$ and $\gamma_{i,i+1}$, $i=1,\dots, 4$, (indexation modulo~4). On~$\mathbb{P}_1(\mathbb{C})$, we have two relations, an exterior one and an interior one, namely,
\begin{gather*}
r_{\rm ext}\colon  \ \gamma_{1,2}\gamma_{2,3}\gamma_{3,4}\gamma_{4,1}={\star}_1 \quad \mbox{(the trivial loop based in $s_1$)},\\
r_{\rm int}\colon \ \gamma_{1,1}\gamma_{1,2}\gamma_{2,2}\gamma_{2,3}\gamma_{3,3}\gamma_{3,4}\gamma_{4,4}\gamma_{4,1}={\star}_1.
\end{gather*}
The local fundamental subgroupoid is generated by the loops~$\gamma_{i,i}$, and is a disjoint union of four monogeneous groups.

\textbf{Representations of the groupoid $\boldsymbol{\pi_1(X,S)}$.}
 A representation of $\pi_1(X,S)$ in a group~$G$ is a~morphism of groupoids $\rho$ from $\pi_1(X,S)$ into~$G$. The group~$G$ is here a groupoid with only one object, whose morphisms are the elements of~$G$. Therefore $\rho$ is characterized by its action on the morphisms of~$\pi_1(X,S)$.

\textbf{Analytic representations of the groupoid $\boldsymbol{\pi_1(X,S)}$ in $\boldsymbol{G}$ induced by a connection~$\boldsymbol{\nabla}$ in~$\boldsymbol{\mathcal{E}}$.}
\begin{enumerate}\itemsep=0pt
\item[--] For each object $s=(p,d)$, we consider a fundamental system of holomorphic solutions~$X_s$ in a neighborhood of $s$ in $X$, i.e., in a small sector at $p$ around the direction $d$, admitting an asymptotic expansion at~$p$.
\item[--]  At each morphism $\gamma_{i,j}$ joining $s_i$ to $s_j$, corresponds a connection matrix~$M_{i,j}$ between the fundamental systems of solutions~$X_i$ and~$X_j$ chosen at~$s_i$ and~$s_j$ def\/ined by
\begin{gather}
X_j=\widetilde{X_i}^{\gamma_{i,j}}M_{i,j},\label{rep_morphism}
\end{gather}
where $\widetilde{X_i}^{\gamma_{i,j}}$ is the analytic continuation of~$X_i$ along~$\gamma_{i,j}$. With this notation
\begin{gather*}
\rho(\gamma_{i,j}\gamma_{j,k})=\rho(\gamma_{i,j})\rho(\gamma_{j,k}).
\end{gather*}
\end{enumerate}

\subsection{The character variety}

\begin{Definition} Let $\rho\colon \pi_1(X,S)\rightarrow G$ and $\rho'\colon \pi_1(X,S)\rightarrow G$ be two analytic representations of $\pi_1(X,S)$.
Let $\rho(\gamma_{i,j})=M_{i,j}$, and $\rho'(\gamma_{i,j})=M'_{i,j}$, where $\gamma_{i,j}$ is a morphism from~$s_i$ to~$s_j$.
The two representations $\rho$ and $\rho'$ are equivalent if and only if for each object~$s_i$, there exists $N_i$ in $G$ such that
\begin{gather*}
M'_{i,j}=N_iM_{i,j}N_j^{-1}.
\end{gather*}
\end{Definition}

Therefore, if we change the choice of the fundamental system attached to each object~$s_i$, we obtain a new equivalent representation. All the representations induced by~$\Delta$ are equivalent. Furthermore, if~$\Delta$ and~$\Delta'$ are gauge equivalent, their representations are equivalent.
The class~$[\rho]$ only depends on~$\nabla$.

\begin{Definition}
Let $\mathcal{R}(S)$ (resp.~$\mathcal{R}(S)^{\rm loc}$) be the space of the analytic representations of $\pi_1(X,S)$ (resp.\ of~$\pi_1^{\rm loc}(X,S)$) induced by some $\Delta$ in $\mathcal{E}$, and $\sim$ the above equivalence relation on $\mathcal{R}(S)$.
\begin{itemize}\itemsep=0pt
  \item The character variety $\chi(S)$ is the quotient $\mathcal{R}(S)/{\sim}$.
  \item In the same way, the local character variety is $\chi(S)^{\rm loc}=\mathcal{R}(S)^{\rm loc}/{\sim}$. The morphism $\pi$ from $\chi(S)$ to $\chi(S)^{\rm loc}$ is induced by restriction of the representations to $\pi_1^{\rm loc}(X,S)$.
\end{itemize}
\end{Definition}

\textbf{Normalized representations.}
We construct a ``good'' representative of $[\rho]$ in $\chi(S)$ by using the following process:
\begin{enumerate}\itemsep=0pt
\item[--] We choose freely a fundamental system of solutions $X_1$ in $s_1$.
\item[--] In $s_2$, we choose $X_2$ to be the analytic continuation of $X_1$ along the path~$\gamma_{1,2}$, then we choose~$X_3$ by analytic continuation of $X_2$ along~$\gamma_{2,3}$, and f\/inally $X_4$ by analytic continuation along $\gamma_{3,4}$. With these choices, we have
\begin{gather*}\rho(\gamma_{1,2})=\rho(\gamma_{2,3})=\rho(\gamma_{3,4})=I,\end{gather*}
and from the exterior relation, we obtain $\rho(\gamma_{1,4})=I$.
The representation~$\rho$ is now characterized by 4 matrices $M_{i,i}=\rho(\gamma_{i,i})$. From the interior relation, we have
\begin{gather*}
M_{1,1}M_{2,2}M_{3,3}M_{4,4}=I.
\end{gather*}
A change in the initial choice of $X_1$ will give rise to 4 matrices related to the previous ones by a common conjugation.
Finally, we have characterized~$[\rho]$ in~$\chi(S)$ by the data of~3 matrices~$M_{i,i}$, $i=1,2,3$ up to a~common conjugation, as in the usual presentation with only one base point. Nevertheless, this groupoid point of view is more convenient, f\/irst for computing the isomonodromic dynamics, but also to get an extension to the irregular cases.
\end{enumerate}

\subsection{The character variety in trace coordinates}
 We will now describe the af\/f\/ine algebraic structure of $\chi(S)$ thanks to the following lemma:

\begin{Lemma}[Fricke lemma~\cite{Mag}]\label{Fric} Given $3$ matrices $M_1$, $M_2$ and $M_3$ in ${\rm SL}(2,\mathbb{C})$ we denote
\begin{gather*}
a_i=\operatorname{tr}(M_i),\qquad x_{i,j}=\operatorname{tr}(M_iM_j),\qquad x_{i,j,k}=\operatorname{tr}(M_iM_jM_k),
\end{gather*}
where the indices are $2$ by $2$ distincts. Since the trace map is invariant under cyclic permutations, we have $3$ coordinates $a_i$, $3$ coordinates~$x_{i,j}$ and $2$ coordinates~$x_{i,j,k}$. We have the following relations
\begin{gather*}%\label{Fricke}
x_{1,2,3}+x_{1,3,2} = a_1x_{2,3}+a_2x_{3,1}+a_3x_{1,2}-a_1a_2a_3:=P, \\
x_{1,2,3}\times x_{1,3,2} = a_1^2+a_2^2+a_3^2+x_{1,2}^2+x_{2,3}^2+x_{3,1}^2+x_{1,2}x_{2,3}x_{3,1}\\
\hphantom{x_{1,2,3}\times x_{1,3,2} =} {}-a_1a_2x_{1,2}-a_2a_3x_{2,3}-a_3a_1x_{3,1}-4:=Q.
\end{gather*}
Therefore, $x_{1,2,3}$ and $x_{1,3,2}$ are the two solutions of the equation $X^2-PX+Q=0$.
\end{Lemma}

Let $M_4=(M_1M_2M_3)^{-1}$ and $a_4=\operatorname{tr}(M_4)$. We have $a_4=x_{1,2,3}$. From the Fricke lemma, we obtain the following relation in $\mathbb{C}^4\times\mathbb{C}^3$
\begin{gather*}F\colon \ a_4^2-Pa_4+Q=0.\end{gather*}
We call it the ``Fricke hypersurface''.
It is a quartic in $\mathbb{C}^7$ endowed with the coordinates $(a_1,a_2,a_3,a_4,x_{1,2},x_{2,3},x_{3,1})$ and, with respect to the last 3 coordinates, a family of cubics~$F_a$ indexed by $a=(a_1,a_2,a_3,a_4)$.

By applying this lemma to the matrices $M_{i,i}=\left(\begin{matrix}\alpha_i&\beta_i\\ \gamma_i&\delta_i\end{matrix}\right)$ for $i=1,2,3$ of a~normalized representation, the trace coordinates
 $(a_1,a_2,a_3,a_4,x_{1,2},x_{2,3},x_{3,1})$ def\/ine a morphism
 \begin{gather*}
 T\colon \ \chi(S)\rightarrow F.
 \end{gather*}
We consider the open set $\chi(S)^*$ def\/ined by the following conditions:
\begin{enumerate}\itemsep=0pt
\item[(i)] each matrix $M_{i,i}$ is semi-simple;
\item[(ii)] one of them (say $M_{1,1}$) is dif\/ferent from~$\pm I$;
\item[(iii)] the two others satisfy $\beta_2\gamma_2\beta_3\gamma_3\neq 0$.
\end{enumerate}
We set $\chi(S)^{*{\rm loc}}=\pi^*\chi(S)^*$.

\begin{Proposition} The morphism $T$ is an isomorphism from~$\chi(S)^*$ onto~$F^*:=(a_1\neq\pm 2)$.
The restriction of $T$ on $\chi(S)^{*{\rm loc}}$ is an isomorphism onto $\mathbb{C}{\setminus}\{\pm 2\}\times \mathbb{C}^3$ and we have~$T\pi=p_1 T$, where~$p_1$ is the first projection $\mathbb{C}^4\times \mathbb{C}^3\rightarrow\mathbb{C}^4$.
\end{Proposition}

\begin{proof}
By using a conjugation, we may suppose that
\begin{gather*}
M_{1,1}=\left(\begin{matrix}\alpha_1&0\\ 0&\alpha_1^{-1}\end{matrix}\right) \quad \mbox{with} \ \ \alpha_1\neq\pm1,\qquad M_{i,i}=\left(\begin{matrix}\alpha_i&\beta_i\\ \gamma_i&\delta_i\end{matrix}\right) \quad \mbox{for} \ \ i=2, 3.
\end{gather*}
This writing is not still unique: we may use a conjugation by a diagonal matrix~$D$. Since the center do not act, we may suppose that $\det(D)=1$, i.e.,
\begin{gather*}
D=D_\alpha=\left(\begin{matrix}\alpha&0\\ 0&\alpha^{-1}\end{matrix}\right).
\end{gather*}

\begin{Lemma}\label{triplets}
Under the conditions $(i)$, $(ii)$ and $(iii)$ defining $\chi(S)^*$, two triples
 $(M_1,M_2,M_3)$ and $(M'_1,M'_2,M'_3)$ are in the same orbit for the action of the group $\{D_\alpha,\, \alpha\in\mathbb{C}^*\}$ if and only if
\begin{gather*}
\alpha_1=\alpha'_1,\qquad \alpha_2=\alpha'_2,\qquad \delta_2=\delta'_2,\qquad \alpha_3=\alpha'_3,\qquad \delta_3=\delta'_3,\\ \beta_2\gamma_3=\beta'_2\gamma'_3,\qquad \gamma_2\beta_3=\gamma'_2\beta'_3.
\end{gather*}
\end{Lemma}

\begin{proof}
Since
\begin{gather*}
D_\alpha M_{i,i}D_\alpha^{-1}=\left(\begin{matrix}\alpha_i&\alpha^2\beta_i\\ \alpha^{-2}\gamma_i&\delta_i\end{matrix}\right),
\end{gather*}
the condition is necessary. Suppose now that this condition holds for two triples. Since $\beta_2\gamma_3\neq 0\neq \beta'_2\gamma'_3$ we can choose $\alpha$ such that $\alpha^2=\frac{\beta'_2}{\beta_2}=\frac{\gamma_3}{\gamma'_3}.$ We also have $\beta_2\gamma_2=1-\alpha_2\delta_2=1-\alpha'_2\delta'_2=\beta'_2\gamma'_2$. Therefore $\alpha^{-2}=\frac{\gamma'_2}{\gamma_2}=\frac{\beta_3}{\beta'_3}$, which proves that $D_\alpha M_{i,i}D_\alpha^{-1}=M'_{i,i}$.
\end{proof}

Now we have to solve in ${\rm SL}_2$ the system
\begin{subequations}\label{s(1)-s(7)}
\begin{gather}
     \alpha_1+\alpha_1^{-1}=a_1, \label{s(1)}\\
     \alpha_2+\delta_2=a_2, \label{s(2)}\\
     \alpha_3+\delta_3=a_3, \label{s(3)}\\
     \alpha_1\alpha_2+\alpha_1^{-1}\delta_2=x_{1,2}, \label{s(4)}\\
     \alpha_1\alpha_3+\alpha_1^{-1}\delta_3=x_{1,3}, \label{s(5)}\\
     \alpha_2\alpha_3+\beta_2\gamma_3+\gamma_2\beta_3+\delta_2\delta_3=x_{2,3}, \label{s(6)}\\
     \alpha_1\alpha_2\alpha_3+\alpha_1\beta_2\gamma_3+\alpha_1^{-1}\gamma_2\beta_3+\alpha_1^{-1}\delta_2\delta_3
     =x_{1,2,3}. \label{s(7)}
  \end{gather}
\end{subequations}
We choose one of the two solutions of equation \eqref{s(1)}. Since $\alpha_1-\alpha_1^{-1}\neq 0$, we obtain from equations~\eqref{s(2)},~\eqref{s(3)},~\eqref{s(4)} and~\eqref{s(5)} a unique solution for $\alpha_2$, $\delta_2$, $\alpha_3$, $\delta_3$. Equations \eqref{s(6)},~\eqref{s(7)} def\/ine a linear system in the 2 variables $(\beta_2\gamma_3,\gamma_2\beta_3)$ of maximal rank if
$\alpha_1-\alpha_1^{-1}\neq 0$. We obtain a unique solution for $\alpha_2$, $\delta_2$, $\alpha_3$, $\delta_3$,
$\beta_2\gamma_3$, $\gamma_2\beta_3$, and therefore a unique triple $(M_1,M_2,M_3)$ up to conjugation according to the preliminary remark.
Note that this triple is not necessarily in~${\rm SL}_2$: the compatibility condition corresponds to the Fricke relation.

Now if we begin with the second solution of~\eqref{s(1)}, the new matrix~$M'_1$ satisf\/ies $M'_1=PM_1P^{-1}$, where~$P$ is the matrix of the transposition. The system~\eqref{s(2)}--\eqref{s(7)} has a unique solution $M'_i$ under the same assumption $\alpha_1-\alpha^{-1}\neq 0$. Since we know that $PM_iP^{-1}$ is another pre-image for $T$, we have: $M'_i=PM_iP^{-1}$. Therefore this second solution is conjugated to the f\/irst one by $P$, and we obtain a unique pre-image of a point in~$F^*$ in~$\chi(S)^*$.
For the second part of the statement, if each matrix~$M_i$ is a semi-simple one, the trace of~$M_i$ characterizes the conjugation class of~$M_i$ in ${\rm SL}_2$.
\end{proof}

\subsection[The dynamics on $\chi(S)$]{The dynamics on $\boldsymbol{\chi(S)}$}

We set
\begin{gather*}
\chi=\cup_{S\in\mathcal{C}}  \chi(S),
\end{gather*}
where $S$ belongs to the space $\mathcal{C}$ of the conf\/igurations of 4 distinct points in the plane. This f\/ibration is endowed with a f\/lat connection (the isomonodromic connection) whose local trivia\-li\-sations are def\/ined by identifying the generators~$\gamma_{i,j}(S)$ and $\gamma_{i,j}(S')$, for~$S'$ near from~$S$. We want to compute the monodromy of this connection on a f\/iber~$\chi(S)$.

The fundamental group $\pi_1(\mathcal{C},[S])$ is the
pure braid group~$P_4$. It is generated by the 3 elements~$b_1$, $b_2$, $b_3$, where~$b_i$ is the pure braid between~$s_i$ and $s_{i+1}$, with the relation $b_1b_2b_3=\operatorname{id}$ (note that the cross ratio induces an isomorphism from~$\mathcal{C}$ on~$P^1(C){\setminus}\{0,1,\infty\}$).

The generators $b_i$ induce an isomorphism from $P_4$ to the mapping class group of the disc punctured by 4 holes, with a base point on their boundaries. This interpretation allows us to construct an action from~$P_4$ on the groupoid~$\pi_1(X,S)$. We denote by~$h_1$, $h_2$, $h_3$ the images of the braids in $\operatorname{Aut}_0(\pi_1(X,S))$.

The automorphisms $h_i$ act on $\mathcal{R}$ by $h_{i*}\colon \rho\mapsto\rho\circ h_i$ and an inner automorphism sends~$\rho$ on an equivalent representation. Therefore each $[h_i]$ in $\operatorname{Out}_0^*(\pi_1(X,S))$ acts on $\chi(S)$.

Looking at the picture of the groupoid, we immediately obtain:
\begin{Proposition}\quad{\samepage
\begin{enumerate}\itemsep=0pt
  \item[$1)$] $h_1(\gamma_{i,i})=\gamma_{i,i}$, $i=1,\dots, 4$;
  \item[$2)$] $h_1(\gamma_{3,2})=\gamma_{3,2}\gamma_{2,1}\gamma_{1,1}\gamma_{1,2}\gamma_{2,2}$;
  \item[$3)$] $h_1(\gamma_{1,2})=\gamma_{1,2}$, $h_1(\gamma_{3,4})=\gamma_{3,4}$.
\end{enumerate}
We have similar expressions for $h_2$ and $h_3$ by cyclic permutations of the indices.}
\end{Proposition}

Now we compute the action~$h_{i*}$ on $\chi(S)$ in three steps. Let~$[\rho]$ in $\chi(S)$ given by a~norma\-li\-zed representation~$\rho$, and therefore by 3 matrices $M_i$ up to a common conjugation.
For each generator~$b_i$,
\begin{enumerate}\itemsep=0pt
  \item[1)] we compute $\rho\circ h_i$ on the generating morphisms $\gamma_{i,j}$. The representation $\rho\circ h_i$ is not yet a normalized one;
  \item[2)] we normalize $\rho\circ h_i$ in a new equivalent representation $\rho'$, by changing the representation of the objects. Let $M'_i$ be the matrices related to $\rho'$;
  \item[3)] we compute $\operatorname{tr}(M'_iM'_j)$ as expressions in the $\operatorname{tr}(M_iM_j)$'s in order to write $h_{i*}$ in the trace coordinates
  $a_i$, $x_{1,2}$, $x_{2,3}$, $x_{3,1}$.
\end{enumerate}
For this last step, we will make use of an extended version of the Fricke lemma:

\begin{Lemma}[extended Fricke lemma]\label{Fricke.bis} We follow the notations of the Fricke Lemma~{\rm \ref{Fric}}. We have
\begin{gather*}
\operatorname{tr}\big(M_1M_2M_1^{-1}M_3\big)= -x_{1,2}x_{1,3}-x_{2,3}+a_1a_4+a_2a_3,\\
\operatorname{tr}\big(M_1M_2M_1M_2^{-1}M_1^{-1}M_3\big)=x_{1,2}^2x_{1,3}\!+x_{1,2}x_{2,3}\!-x_{1,3}\!-x_{1,2}(a_1a_4+a_2a_3)\!+(a_1a_3+a_2a_4).
\end{gather*}
\end{Lemma}

\begin{proof}
We only make use of the relation $\operatorname{tr}(AB)+\operatorname{tr}(AB^{-1})=\operatorname{tr}(A)\operatorname{tr}(B)$:
\begin{gather*}
\operatorname{tr}\big(M_1M_2M_1^{-1}M_3\big)=\operatorname{tr}\big(M_3M_1M_2M_1^{-1}\big)
 = a_4a_1-\operatorname{tr}(M_3M_1M_2M_1)\\
\qquad{}= a_4a_1-\big(x_{3,1}x_{1,2}-\operatorname{tr}\big(M_3M_2^{-1}\big)\big)
 = a_4a_1-(x_{3,1}x_{1,2}-(a_2a_3-x_{3,2}))\\
\qquad{}
= -x_{1,2}x_{1,3}-x_{2,3}+a_1a_4+a_2a_3,
\\
\operatorname{tr}\big(M_1M_2M_1M_2^{-1}M_1^{-1}M_3\big)=\operatorname{tr}\big(M_2^{-1}M_1^{-1}M_3M_1\cdot M_2M_1\big)\\
\qquad {}= \operatorname{tr}\big(M_2^{-1}M_1^{-1}M_3\cdot M_1\big)x_{1,2}-\operatorname{tr}\big(M_2^{-1}M_1^{-1}M_3M_2^{-1}\big)\\
\qquad {}= \operatorname{tr}(M_2^{-1}M_1^{-1}M_3)a_1x_{1,2}-\operatorname{tr}(M_2^{-1}M_1^{-1}\cdot M_3M_1^{-1})x_{1,2}\\
\qquad \quad{} -\operatorname{tr}\big(M_2^{-1}M_1^{-1}M_3\big)a_2+\operatorname{tr}\big(M_2^{-1}M_1^{-1}M_3M_2\big)\\
\qquad {}= (x_{1,2}a_3-a_4)a_1x_{1,2}-x_{1,2}^2(a_1a_3-x_{1,3})+x_{1,2}x_{2,3}-(x_{1,2}a_3-a_4)a_2+(a_1a_3-x_{1,3})\\
\qquad {}= x_{1,2}^2x_{1,3}+x_{1,2}x_{2,3}-x_{1,3}-x_{1,2}(a_1a_4+a_2a_3)+(a_1a_3+a_2a_4).\tag*{\qed}
\end{gather*}\renewcommand{\qed}{}
\end{proof}

\begin{Proposition} Let $x'_{i,j}=\operatorname{tr}(M'_iM'_j).$ In these trace coordinates, $h_{1*}$ is given by
\begin{gather*}
x'_{1,2} = x_{1,2}, \\
x'_{2,3} = -x_{1,2}x_{1,3}-x_{2,3}+(a_1a_4+a_2a_3),\\
x'_{1,3} = x_{1,2}^2x_{1,3}+x_{1,2}x_{2,3}-x_{1,3}-x_{1,2}(a_1a_4+a_2a_3)+(a_1a_3+a_2a_4).
\end{gather*}
We obtain~$h_{2*}$ and~$h_{3*}$ by a cyclic permutation indices~$(+1)$ of the indices~$1$,~$2$ and~$3$.
\end{Proposition}

\begin{proof}
We have
\begin{gather*}
\rho\circ h_1(\gamma_{1,2})=\rho(\gamma_{1,2})=I,\\
\rho\circ h_1(\gamma_{3,4})=\rho(\gamma_{3,4})=I,\\
\rho\circ h_1(\gamma_{3,2})=\rho(\gamma_{3,2}\gamma_{2,1}\gamma_{1,1}\gamma_{1,2}\gamma_{2,2})=M_1M_2.
\end{gather*}
Therefore $X_3=\widetilde{X_2}^{h_1(\gamma_{2,3})}(M_1M_2)^{-1}$, and we normalize $\rho\circ h_1$ by setting: $X'_1=X_1$, $X'_2=X_2$ and $X'_3=X_3\cdot M_1M_2$ in order to obtain a representation $\rho'$ equivalent to~$\rho\circ b_1$ which satisf\/ies $\rho'(\gamma_{2,3})=I$. This representation $\rho'$ is characterized by the 3 matrices:
\begin{gather*}M'_1=M_1,\qquad M'_2=M_2,\qquad M'_3=(M_1M_2)^{-1}M_3(M_1M_2).
\end{gather*}
The statement of the proposition is obtained from the extended Fricke Lemma~\ref{Fricke.bis}.
\end{proof}

By this way, we reach the same expressions of the dynamics as S.~Cantat and F.~Loray in~\cite{CanLo}.
The study of this dynamics allows them to give a new proof of the irreducibility of the Painlev\'e~VI equation.
This is an important motivation to extend the description of this dynamics to the non regular cases.

\section[An irregular example (towards $P_{\rm V}$)]{An irregular example (towards $\boldsymbol{P_{\rm V}}$)}

\subsection{Standard facts about irregular singularities}
      We consider a rank $n$ linear dif\/ferential system at $z=0$
      \begin{gather*}
      \Delta\colon \ z^{r+1}\frac{dY}{dz}=A(z)\cdot Y, \qquad  A \ \mbox{holomorphic at} \ 0,
      \end{gather*}
      where $A(z)$ takes its values in $\mathfrak{gl}(n,\mathbb{C})$.\footnote{The theory can be extended to any complex reductive Lie algebra, up to some technical complications, see~\cite{BaVa}.}
      The integer~$r$ is positive, and equal to~0 for a~Fuchsian system. We suppose here that $r>0$, and that the eigenvalues of $A_0=A(0)$ are non vanishing distincts complex numbers. We f\/ix here $\Lambda_0=\operatorname{diag} (\lambda_i)$ a normal form of~$A_0$ in the Cartan subalgebra $\mathcal{T}_0$ of the diagonal matrices, i.e., we choose an ordering of its eigenvalues.
The formal local meromorphic classif\/ication is given by

\begin{Proposition}\label{formal_sol}\quad
\begin{enumerate}\itemsep=0pt
\item[$1.$] Up to a local ramified formal meromorphic gauge equivalence, we have
\begin{gather*}
\Delta\sim_0\frac{dX}{dt}=\left(\frac{dQ}{dt}+\frac{L}{t}\right)\cdot X,
\end{gather*}
where $z=t^\nu$, $Q$ $($the ``irregular type''$)$ $=\frac{\Lambda_0}{t^r}+\cdots +\frac{\Lambda_{r-1}}{t}$, and the matrices $\Lambda_i$ and $L$ $($the residue matrix$)$ are diagonal matrices.
For a fixed $\Lambda_0$, the pair $(Q,L)$ is unique in $\mathcal{T}_0\times\mathcal{T}_0/\mathcal{T}_0(\mathbb{Z})$.
\item[$2.$] Let $F_0$ be a conjugation between $A_0$ and $\Lambda_0$: $A_0=F_0\Lambda_0F_0^{-1}$. The system $\Delta$ has a~formal fundamental solution
\begin{gather*}\widehat{X}=\widehat{F}(t)t^{L}\exp Q\qquad \mbox{with} \quad \widehat{F}(0)=F_0.
\end{gather*}
\end{enumerate}
\end{Proposition}

For a f\/ixed $\Lambda_0$, there are already two ambiguities in the above writing of $\widehat{X}$:
\begin{enumerate}\itemsep=0pt
\item[--] the choice of $F_0$: we may change $F_0$ with $F_0D$, and therefore $\widehat{X}$ with $\widehat{X}D$, where $D$ belongs to the centralizer of $\Lambda_0$;

\item[--] the choice of a branch for the argument, and hence for $\log t$ and $t^{L}$.
\end{enumerate}

We suppose now that we are in the unramif\/ied case: $\nu=1$.
\begin{Definition}\quad
\begin{enumerate}\itemsep=0pt
  \item A \emph{separating ray} is a ray $\operatorname{arg}(z)=\tau$ such that there exists a pair of eigenvalues $(\lambda_j,\lambda_k)$ of~$\Lambda_0$ satisfying:
      $z^r(\lambda_j-\lambda_k)\in i\mathbb{R}^+$ for $\operatorname{arg}(z)=\tau$.
  \item A \emph{singular ray} is a ray $\operatorname{arg}(z)=\sigma$ such that there exists a pair of eigenvalues $(\lambda_j,\lambda_k)$ satisfying:
      $z^r(\lambda_j-\lambda_k)\in \mathbb{R}^-$ for $\operatorname{arg}(z)=\sigma$.\footnote{Since in many references, the def\/initions of Stokes and anti-Stokes rays are exchanged, we do not use this terminology here.}
  \item A \emph{regular sector} is an open sector $\mathcal{S}_d$ of angle $\pi/r$ bisected by $d$ such that~$d$ is not a~singular ray (or equivalently, such that its edges are not separating rays).
\end{enumerate}
\end{Definition}

We can remark that:
\begin{itemize}\itemsep=0pt
  \item If $\operatorname{arg}(z)=\tau$ is a separating (resp.\ singular) ray then its opposite is also a separating (resp.\ singular) ray.
  \item A non singular ray is a ray on which the formal solution admits a unique sum, for the summation theory.\footnote{We only need here $k$-summation theory, with $k=r$.}
  \item A separating ray is a ray on which the asymptotic of a general solution (a~linear combination of the columns of $\widehat{X}$) changes.
      \item The knowledge of the $\mu$ separating (resp.\ singular) rays in a regular sector generates the complete knowledge of all the separating (resp.\ singular) rays, by considering their opposites, and the ramif\/ication by~$z^r$.
  \item The generic case is the situation in which there exists exactly one pair of eigenvalues $(\lambda_j(\nu),\lambda_k(\nu))$ def\/ining each separating ray~$\tau_\nu$. In this case, we have $\mu=n(n-1)/2$ separating rays in a regular sector, and $m=n(n-1)r=2r\mu$ separating (resp.\ singular) rays in~$S^1$.
\end{itemize}

\begin{Theorem}\label{sol_regulieres} On a regular sector $S_d$ containing the $\mu$ separating rays $\tau_\nu,\dots,\tau_{\nu+\mu-1}$, there exists a unique holomorphic fundamental system of solutions
        $X_d$ admitting the asymptotic expansion~$\widehat{X}$. Furthermore $X_d$ can be extended to a solution $($with the same asymptotic$)$ on the sector~$S_\nu$ delimited by the two nearest separating rays~$\tau_{\nu-1}$ and~$\tau_{\nu+\mu}$ outside~$S_d$.
\end{Theorem}
There exists two proofs of this fact using either the asymptotic theory (see~\cite{BJL}), or the summation theory (see~\cite{MR2}).

For $\nu=1$ to~$m$, Theorem~\ref{sol_regulieres} gives us a unique solution~$X_\nu$ on the large sector delimited by~$\tau_{\nu-1}$ and $\tau_{\nu+\mu}$ admitting~$\widehat{X}$ as asymptotic expansion on~$\mathcal{S}_\nu$.
\begin{Definition}
The $m$ Stokes multipliers $U_\nu$ are def\/ined on $\mathcal{S}_\nu\cap \mathcal{S}_{\nu+1}$ by
\begin{gather*}
X_{\nu}=X_{\nu+1}\cdot U_\nu,
\end{gather*}
with a $m$-periodic indexation.
\end{Definition}

In the generic case (the support of each singular ray reduces to a unique pair of eigenvalues), the constant matrices $U_\nu$, $\nu=1,\dots, m-1$ are transvection matrices: the diagonal entries are equal to~1, and the unique non vanishing coef\/f\/icient of\/f the diagonal is the coef\/f\/icient in position $(j(\nu),k(\nu))$
where the separating ray~$\tau_\nu$ (the only ray in $\mathcal{S}_\nu{\setminus} \mathcal{S}_\nu\cap \mathcal{S}_{\nu+1}$) is def\/ined by the pair $(\lambda_j,\lambda_k)$.
The diagonal of the matrix~$U_m$ is $\exp(2i\pi L)$, where~$L$ is the residue coef\/f\/icient of~$A(x)$ after diagonalisation.

\subsection[The wild fundamental groupoid for the class of connections $(0,0,1)$]{The wild fundamental groupoid for the class of connections $\boldsymbol{(0,0,1)}$}

We consider a meromorphic $\mathfrak{sl}_2(\mathbb{C})$-connection $\nabla$ on ${\mathbb P}_1(\mathbb{C})$, admitting $2$ regular singular points (say~0 and~1) and an irregular one at $\infty$. This corresponds to an element of the family indexed by~$(0,0,1)$ in the classif\/ication of M.~van der Put and K.~Saito in~\cite{VdPSai}.
Locally in a coordin\-a\-te~$z$ centered at $\infty$, $\nabla$ is given by the following system
\begin{gather*}
\Delta\colon \ \frac{dX}{dz}=z^{-2}\sum_{i\geq 0}A_iz^i\cdot X.
\end{gather*}
The initial part $A_0\neq 0$ is a non trivial semi-simple element of $\mathcal{G}={\mathfrak{sl}}_2$.
We choose a normal form $\Lambda_0=\operatorname{diag} (\lambda_0,-\lambda_0)$, $\lambda_0\neq 0$, in $\mathcal{T}_{0}$ (the Cartan algebra of diagonal matrices) in the conjucacy class of $A_0$, i.e., an ordering of its eigenvalues. Its centralizer in ${\rm SL}_2(\mathbb{C})$ is
\begin{gather*}
C(\Lambda_0)=C(\mathcal{T}_0)=\left\{D_\alpha= \left(\begin{matrix}\alpha&0\\ 0&\alpha^{-1}\end{matrix}\right),\, \alpha\in \mathbb{C}^* \right\}.
\end{gather*}
We have here two singular rays $\sigma_1(\Lambda_0)$ and $\sigma_2(\Lambda_0)=-\sigma_1(\Lambda_0)$, and two separating rays $\tau_1(\Lambda_0)$ and $\tau_2(\Lambda_0)=-\tau_1(\Lambda_0)$.

\begin{figure}[t]\centering
\includegraphics[scale=0.45]{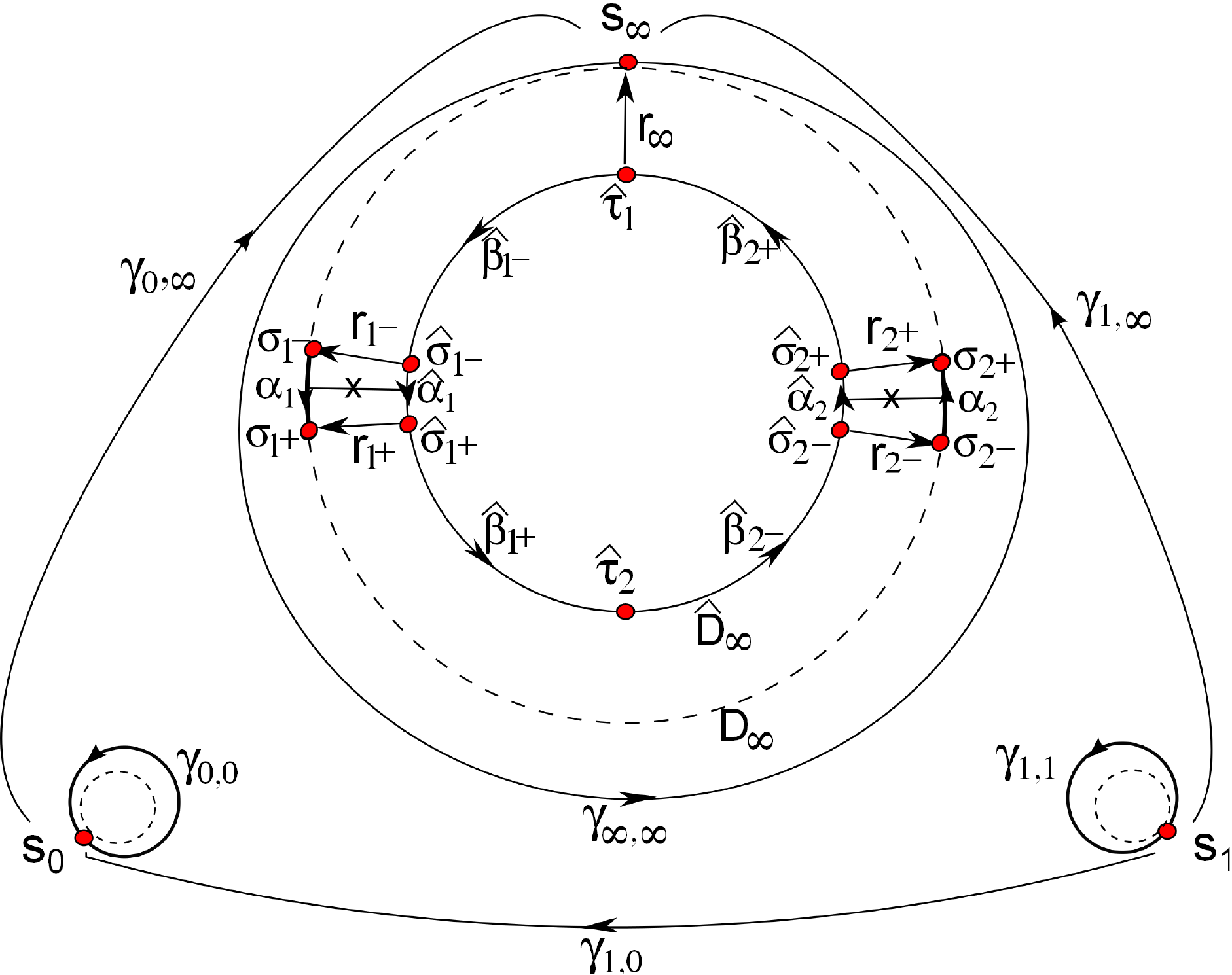}
\caption{The groupoid $\pi_1(X,S)$ (irregular case).}\label{groupoid-PV}
\end{figure}

In order to construct the wild fundamental groupoid (for a f\/ixed $\Lambda_0$), we f\/irst use a real blowing up at each singularity in $\mathbb{P}^1$ and we obtain a variety $X$ with 3 exceptional divi\-sors~$D_0$,~$D_1$ and~$D_\infty$ (circles in dotted lines in Fig.~\ref{groupoid-PV}). As in the previous classical case, we choose a base point~$s_0$,~$s_1$ and $s_\infty$ on each of them, and we consider the morphism (path up to homotopy)~$\gamma_{i,j}$ joining~$s_i$ to~$s_j$. The paths~$\gamma_{i,i}$ are homotopic in~$X$ to the curves~$D_i$. We choose~$s_\infty$ such that it corresponds to a non singular ray.

Since we also have to consider the continuation of $\widehat{X}$ (the formal monodromy along an arc is induced by the substitution $z\mapsto ze^{i\theta}$), we introduce a second copy $\widehat{D}_\infty$ of $D_\infty$ inside the f\/irst one, with a base point $\widehat{\tau}_1$ which is the separating ray between~$\sigma_2$ and~ $\sigma_1$.
 For each singular direction $\sigma_i$ (denoted in the picture below by a ray with a cross in the annulus between $\widehat{D}_\infty$ and $D_\infty$), we add two loops delimited by two rays~$r_i^-$ and~$r_i^+$ which are non singular and non separating, and two arcs~$\alpha_i$ on $D_\infty$ and $\widehat{\alpha}_i$  and $\widehat{D}_\infty$ of opening strictly lower than~$\pi$, bisected by the singular rays. Let $\widehat{\sigma}_i^-$ and~$\sigma_i^-$ the two points on $\widehat{D}_\infty$ and $D_\infty$ joined by~$r_i^-$, and~$\widehat{\sigma}_i^+$ and~$\sigma_i^+$ joined by~$r_i^+$.
 Finally we put a ray~$r_\infty$ from~$\widehat{\tau}_1$ to~$s_\infty$.

\begin{Remark}
 In the picture above, we put arbitrarily the two base points~$s_\infty$ on~$D_\infty$ and~$\widehat{\tau}_1$ on~$\widehat{D}_\infty$ in the same direction: it is our initial conf\/iguration. Nevertheless, in the dynamical study of the next section, $s_\infty$ will remain f\/ixed, while the separating ray~$\widehat{\tau}_1$ (and all the other data related to~$\Lambda_0$) will move on~$\widehat{D}_\infty$.
 \end{Remark}

\begin{Definition}
The wild fundamental groupoid $\pi_1(X,S(\Lambda_0))$ is the groupoid def\/ined by
\begin{itemize}\itemsep=0pt
  \item the objects $S(\Lambda_0)$: the three points $s_0$, $s_1$ and $s_\infty$, the points $\widehat{\tau}_i$ (separating rays), $\widehat{\sigma}_{i^\pm}$~on~$\widehat{D}_\infty$ around the singular rays $\sigma_i$ (denoted on the f\/igure by a ray with a cross), and the corresponding points $\sigma_{i^\pm}$ on $D_\infty$.
  \item the morphisms: they are generated by the paths $\gamma_{i,j}$ (up to homotopy) between $s_0$, $s_1$, and $s_\infty$ in $X$, the rays $r_\infty$, $r_i^\pm$, the arcs $\alpha_i$ on $D_\infty$, and all the arcs on $\widehat{D}_\infty$: $\widehat{\alpha}_i$ from $\widehat{\sigma}_{i^-}$ to $\widehat{\sigma}_{i^+}$, and the connecting arcs $\widehat{\beta}_{i^\pm}$ as indicated on the f\/igure.
\end{itemize}
The subgroupoid $\pi_1^{\rm loc}(X,S(\Lambda_0))$  is generated by the morphisms~$\gamma_{0,0}$, $\gamma_{1,1}$ and~$\widehat{\gamma}_{1,1}$ (the formal loop based in~$\widehat{\tau}_1$).
\end{Definition}

 We still have two relations $r_{\rm int}$ and $r_{\rm ext}$ between the generating morphisms:
\begin{gather*}
r_{\rm int}\colon \ \gamma_{0,0}\gamma_{0,\infty}\gamma_{\infty,\infty}\gamma_{\infty,1}\gamma_{1,1}\gamma_{1,0}
=\star_0,\\
r_{\rm ext}\colon \ \gamma_{0,\infty}\gamma_{\infty,1}\gamma_{1,0}=\star_0
\end{gather*}
and a new one (the wild relation):
\begin{gather*}
r_{\rm wild}\colon \ \gamma_{\infty,\infty}=(r_\infty)^{-1}
\big(\widehat{\beta}_{1^-}r_{1-}\alpha_1 (r_{1^+})^{-1}
\widehat{\beta}_{1^+}\widehat{\beta}_{2^-}r_{2^-}\alpha_2r_{2^+}\widehat{\beta}_{2^+}\big)r_\infty.
\end{gather*}
The \emph{Stokes loops} (based in $\widehat{\sigma}_{i^-}$) are  $\operatorname{st}_i:=r_{i^-}\alpha_i(r_{i^+})^{-1}\widehat{\alpha}_i^{-1}$, $i=1,2$.
The \emph{formal loop} (based in~$\widehat{\tau}_1$) is   $\widehat{\gamma}_{1,1}=\widehat{\beta}_{1^-}\widehat{\alpha}_1\widehat{\beta}_{1^+}
\widehat{\beta}_{2^-}\widehat{\alpha}_2\widehat{\beta}_{2^+}$.

A representation $\rho$ of this groupoid induced by the dif\/ferential system~$(\Delta)$ is def\/ined in the following way. We f\/irst choose a ``compatible'' representation of the objects:
\begin{itemize}\itemsep=0pt
  \item We choose analytic fundamental systems $X(s_0)$, $X(s_1)$, at $s_0$, $s_1$, as in the regular case (i.e., we choose a~logarithmic branch in the corresponding direction);
  \item We choose a formal fundamental system $\widehat{X}_\infty$ given by Proposition~\ref{formal_sol}. For each object~$\widehat{\sigma}_{i^\pm}$, $\widehat{\tau}_i$ on~$\widehat{D}_\infty$, we choose a~formal fundamental system~$X(\widehat{\sigma}_i{^\pm})$, $X(\widehat{\tau}_i)$ by choosing a~determination in the corresponding direction of the formal fundamental solution $\widehat{X}_\infty$.
   \item For each object $\sigma_{i^\pm}$ on~$D_\infty$, we choose an actual sectorial solutions~$X(\sigma_{i^\pm})$, given by Theorem~\ref{sol_regulieres},
   whose asymptotic expansion is some determinacy of the \emph{same formal fundamental solution} $\widehat{X}_\infty$  (this is the compatibility condition).
\end{itemize}
Then, we construct the representations of the generating morphisms in the following way:
\begin{itemize}\itemsep=0pt
  \item We use analytic continuation to represent the morphisms between~$s_0$,~$s_1$, and $s_\infty$, as in the singular regular case (see~(\ref{rep_morphism})).
  \item In the same way, we use the analytic continuation of the formal solutions to represent the morphisms $\widehat{\beta}_{i^\pm}$, and $\widehat{\alpha}_{i}$ on $\widehat{D}_\infty$ between the formal objects: this formal monodromy is def\/ined by the substitution $z\rightarrow ze^{i\theta}$ in the formal expressions.
  \item We also use analytic continuation \emph{preserving the same asymptotic}, to represent the arcs~$\alpha_i$ on~$D_\infty$. Note that the regular sectors given by Theorem~\ref{sol_regulieres} centered on~$r_{i^-}$ and on~$r_{i^+}$ allow us to def\/ine this continuation along~$\alpha_i$
      and~$\alpha_i^{-1}$:  indeed the intersection of these two sectors is a sector of opening $\pi$ delimited by the 2 separating rays, and therefore contains~$\alpha_i$.
  \item We represent $r_{i^\pm}$ and $r_\infty$ by using Theorem~\ref{sol_regulieres}: starting from the representation $X(\widehat{\sigma}_{i^\pm})$ of the formal objects $\widehat{\sigma}_{i^\pm}$, this theorem gives us an actual solution in this direction, which is denoted by $\widetilde{X(\widehat{\sigma}_{i^\pm})}$. The comparison with the representation~$X(\sigma_{i^\pm})$ of the f\/inal object~$\sigma_{i^\pm}$ def\/ines~$\rho(r_{i^\pm})$:
      \begin{gather*}
      X(\sigma_{i^\pm})=\widetilde{X(\widehat{\sigma}_{i^\pm})}\cdot\rho(r_{i^\pm}).
      \end{gather*}
      The representation of the inverse paths~$(r_{i^\pm})^{-1}$ are obtained in the following way: starting from the representation $X(\sigma_{i^\pm})$ of $\sigma_{i^\pm}$ on~$D_\infty$ we use its asymptotic expansion~$\widehat{X(\sigma_{i^\pm})}$ and compare it with the representation $X(\widehat{\sigma}_{i^\pm})$ of the f\/inal object $\widehat{\sigma}_{i^\pm}$ in order to def\/ine
      $\rho((r_{i^\pm})^{-1})$.
\end{itemize}

\begin{Definition}
Two representations of $\pi_1(X,S(\Lambda_0))$ are equivalent if they are obtained by dif\/ferent compatible representations of the objects.
The wild character variety~$\chi(\Lambda_0)$ is the set of the representations of~$\pi_1(X,S(\Lambda_0))$ up to the above equivalence relation.
The local wild character variety $\chi^{\rm loc}(\Lambda_0)$ is the set of the representations of~$\pi_1^{\rm loc}(X,S(\Lambda_0))$ up to the above equivalence relation, and we have a~natural f\/ibration
\begin{gather*}
\pi\colon \  \chi(\Lambda_0)\longrightarrow\chi^{\rm loc}(\Lambda_0).
\end{gather*}
\end{Definition}

\subsection{The normalized representations}

We f\/irst choose a representation $X(\widehat{\tau}_1)$ of  $\widehat{\tau}_1$ (our initial object). Now, we can f\/ix the representations of~$s_0$,~$s_1$ and~$s_\infty$ in a unique way such that
\begin{gather*}
\rho(r_\infty)=\rho(\gamma_{\infty,0})=\rho(\gamma_{\infty,1})=I\quad (=\rho(\gamma_{0,1}) \ \mbox{from the exterior relation}).
\end{gather*}
Then we choose $X(\widehat{\tau}_i)$, $X(\widehat{\sigma}_{i^\pm})$ and $X(\sigma_{i^\pm})$ such that
\begin{gather*}
\rho(\widehat{\beta}_{1^\pm})=\rho(r_{i^\pm})=\rho(\widehat{\alpha}_i)=\rho(\widehat{\beta}_{2^-})=I.
\end{gather*}
There remains f\/ive matrices
\begin{gather*}
M_0=\rho(\gamma_{0,0}), \qquad M_1=\rho(\gamma_{1,1}), \qquad \widehat{M}=\rho(\widehat{\beta}_{2^+}), \qquad
U_1=\rho(\alpha_1), \qquad
U_2=\rho(\alpha_2).
\end{gather*}

\begin{Remark}\quad
\begin{itemize}\itemsep=0pt
\item For such a normalized representation, we also have $U_i=\rho(\operatorname{st}_i)$, $U_i$ is the representation of the Stokes loops. From the above def\/inition of the representation of the paths~$r_i^\pm$, the matrices~$U_i$ are the Stokes multipliers introduced in the previous section. In particular, they are unipotent matrices
    \begin{gather*}
    U_1= \left(\begin{matrix}
1&u_1\\ 0&1\end{matrix}\right) ,\qquad U_2= \left(\begin{matrix}
1&0\\ u_2&1\end{matrix}\right) .
\end{gather*}
\item We also have $\widehat{M}=\rho(\widehat{\gamma}_{1,1})$. Therefore this matrix is a representation of the formal loop. It is a diagonal matrix
\begin{gather*}\widehat{M}= \left(\begin{matrix}\lambda&0\\0&\lambda^{-1}\end{matrix}\right) .\end{gather*}
\item Let $M_\infty:=\rho(\gamma_{\infty,\infty})$.
From the interior relation~$r_{\rm int}$ in the groupoid we have $M_0M_\infty M_1=I$ and from the wild relation $r_{\rm wild}$, $M_\infty=U_1U_2\widehat{M}$. Therefore we have
    \begin{gather*}
    M_0U_1U_2\widehat{M}M_1=I
    \end{gather*}
    and $\rho$ is given by a 4-uple of independent matrices $(M_0,U_1,U_2,\widehat{M})$, where the $U_i$'s are upper and lower unipotent matrices and $\widehat{M}$ is a diagonal matrix.
    \item If we change the choice of the representation of the initial object $X(\widehat{\tau}_1)$ setting
    $X'(\widehat{\tau}_1)=X(\widehat{\tau}_1)\cdot D_\alpha$, $D_\alpha$ in~$C(\mathcal{T}_0)$, the 4-uple
    $(M_0,U_1,U_2,\widehat{M})$ changes by the common conjugacy with~$D_\alpha$.
\end{itemize}
\end{Remark}

Therefore, for a given~$\Lambda_0$, a representation $\rho$ is characterized by a 4-uple $(M_0,U_1,U_2,\widehat{M})$ up to the conjugation by~$C(\mathcal{T}_0)$.
According to the previous description, the character variety is
\begin{gather*}
\chi(\Lambda_0)=\big\{(M_0,U_1,U_2,\widehat{M})_\sim\big\}.
\end{gather*}
Its dimension is $(3+1+1+1)-1=5$. The character variety of the local datas is
\begin{gather*}
\chi^{\rm loc}(\Lambda_0)=\big\{([M_0]_\sim, [M_1]_\sim,\widehat{M})\big\},
\end{gather*}
where $[M_0]_\sim$ is the conjugation class of $M_0$, $[M_1]_\sim$ is the (independent) conjugation class of $M_1$, and $\widehat{M}$ is diagonal.
If $M_0$ and $M_1$ are semi-simple matrices, it is a 3-dimensional variety, and the f\/iber of $\chi(\Lambda_0)\rightarrow\chi^{\rm loc}(\Lambda_0)$ is a 2-dimensional variety.

\textbf{Changing the choice of $\boldsymbol{\Lambda_0}$.} Let $W:=\{\operatorname{id},w\}$ be the group of permutations of two objects. $W$~is isomorphic to the quotient of the subgroup $\{I,P_w\}$ of ${\rm SL}_2$ by $\{\pm I\}$, where
\begin{gather*}
P_w=\left(\begin{matrix}
0&-1\\ 1&0\end{matrix}\right).
\end{gather*}
The conjugation by $P_w$ on ${\mathfrak{sl}}_2$ only depends on the class of $P_w$ in the quotient. Therefore, we denote $c_w(M)=P_w^{-1}MP_w$. Let $w\cdot \Lambda_0:=c_w(\Lambda_0)=-\Lambda_0$.

We consider the fundamental groupoid $\pi_1(X,S(w\cdot\Lambda_0))$ obtained by a new indexation of the singular rays. The objects are
\begin{gather*}
S(w\cdot\Lambda_0)=\big\{s_0,s_1,s_\infty,\sigma_{w(i)^\pm},\widehat{\sigma}_{w(i)^\pm},\widehat{\tau}_{w(i)}\big\}.
\end{gather*}
The \looseness=-1  generating morphisms are also re-indexed according to the new indexation of their origin and end-point.
We obtain an isomorphism of groupoid $\Phi_w$ from $\pi_1(X,S(\Lambda_0))$ to $\pi_1(X,S(w\cdot\Lambda_0))$. Now this new choice of $\Lambda_0$ also modify the choice of the initial representation of the object by $X'(\widehat{\tau}_1)=X(\widehat{\tau}_1)\cdot P_w$, since we change $F_0$ with $F_0P_w$ (see Proposition~\ref{formal_sol}). Therefore the new representation $\rho'$
is obtained from $\rho$ by the conjugation~$c_w$ by~$P_w$. Finally we have an isomorphism $\psi_w$  from $\chi(\Lambda_0)$ to $\chi(w\cdot\Lambda_0)$ which send~$\rho$ on~$\rho'$ def\/ined by the commutative diagram:
\begin{gather*}
\begin{array}{@{}ccc}
  \pi_1(X,S(\Lambda_0) & \stackrel{\rho}{\longrightarrow} & {\rm SL}_2(\mathbb{C}) \\
  \Phi_w\downarrow &  & \downarrow c_w \\
 \pi_1(X,S(w\cdot\Lambda_0) & \stackrel{\rho'}{\longrightarrow} & {\rm SL}_2( \mathbb{C})
\end{array}
\end{gather*}

\begin{Remark}\label{rho'}
Notice that $\rho'$ is characterized by
\begin{gather*}
\big(M_0',U_1',U_2',\widehat{M}'\big)=\big(P_w^{-1}M_0P_w,P_w^{-1}U_{w(1)}P_w,P_w^{-1}U_{w(2)}P_w,P_w^{-1}\widehat{M}P_w\big)\\
\hphantom{\big(M_0',U_1',U_2',\widehat{M}'\big)} {}=\big(\big(M_0^{-1}\big)^t,U_2^t,U_1^t,\widehat{M}^{-1}\big),
\end{gather*}
where $M^t$ denotes the transposed matrix.
\end{Remark}

\subsection[The wild character variety $\chi(\Lambda_0)$ in trace coordinates]{The wild character variety $\boldsymbol{\chi(\Lambda_0)}$ in trace coordinates}

Let suppose that, with our choice of $\Lambda_0$, $U_1$ is an upper unipotent matrix, and $U_2$ a lower one. We have
\begin{gather*}
M_0=\left(\begin{matrix}a_0&b_0\\c_0&d_0\end{matrix}\right),\qquad
U_1=\left(\begin{matrix}
1&u_1\\ 0&1\end{matrix}\right),\qquad U_2= \left(\begin{matrix}
1&0\\ u_2&1\end{matrix}\right),\qquad \widehat{M}= \left(\begin{matrix}\lambda&0\\0&\lambda^{-1}\end{matrix}\right).
\end{gather*}
The action of a diagonal matrix $D_\alpha$ in $C$ is given by
\begin{gather*}
D_\alpha M_0D_\alpha^{-1} = \left(\begin{matrix}a_0&\alpha^2b_0\\
\alpha^{-2}c_0&d_0\end{matrix}\right),\qquad
D_\alpha\widehat{M}D_\alpha^{-1} =\widehat{M},\\
D_\alpha U_1 D_\alpha^{-1} = \left(\begin{matrix}
1&\alpha^2 u_1\\ 0&1\end{matrix}\right),\qquad
D_\alpha U_2 D_\alpha^{-1} = \left(\begin{matrix}
1&0\\ \alpha^{-2}u_2&1\end{matrix}\right).
\end{gather*}

\begin{Lemma}\label{equivalence}
Two datas $(M_0,U_1,U_2,\widehat{M})$ and $(M_0',U_1',U_2',\widehat{M}')$ such that $u_1u_2\neq 0\neq u'_1u'_2$, are equivalent up to $C(\mathcal{T}_0)$ if and only if
\begin{gather*}
\lambda=\lambda',\qquad u_1u_2=u'_1u_2',\qquad u_1c_0=u'_1c'_0,\qquad u_2b_0=u_2'b'_0,\qquad a_0=a_0', \qquad d_0=d'_0.
\end{gather*}
\end{Lemma}
\begin{proof}
Clearly, these quantities are invariant. Suppose now that they are equal. We choose $\alpha$ such that
\begin{gather*}
\alpha^2=\frac{u_1'}{u_1}=\frac{u_2}{u_2'}\quad (\neq 0).
\end{gather*}
From the other relations we obtain: $c_0=\alpha^2c'_0$, $b'_0=\alpha^2b_0$ and f\/inally
$\widehat{M'}=\widehat{M}$, $U_i'=D_\alpha U_i D_\alpha^{-1}$, and
 $M_0'=D_\alpha M_0D_\alpha^{-1}$.
\end{proof}

We consider the 6 coordinates
\begin{gather*}
\lambda\ (\mbox{f\/irst eigenvalue of }\widehat{M}),\qquad t_0=\operatorname{tr}(M_0), \qquad t_1=\operatorname{tr}(M_1),\\
s=\operatorname{tr}(U_1U_2),\qquad x=\operatorname{tr}(M_0U_1U_2),\qquad y=\operatorname{tr}(M_0\widehat{M}).
\end{gather*}
They are invariant by the action of the centralizer $C(\Lambda_0)$ and therefore they induce a map
$T\colon  \chi(\Lambda_0) \rightarrow \mathbb{C}^6$.
The Fricke lemma applied on the triple $(M_0,U_1U_2,\widehat{M})$ def\/ines a~codimension~1 Fricke variety~$F$ given by
$t_1^2-Pt_1+Q = 0$,
with
\begin{gather*}
P=t_0\big(\lambda^{-1}-\lambda+\lambda s\big)+sy+\big(\lambda+\lambda^{-1}\big)x,\\
Q=t_0^2+s^2+\big(\lambda+\lambda^{-1}\big)^2+x^2+\big(\lambda^{-1}-\lambda+\lambda s\big)^2+y^2+xy\big(\lambda^{-1}-\lambda+\lambda s\big).
\end{gather*}
i.e.,
\begin{gather*}
\lambda xys+x^2+y^2+\big(1+\lambda^2\big)s^2-\big(\lambda-\lambda^{-1}\big)xy-t_1sy-t_1\big(\lambda+\lambda^{-1}\big)x-\big(\lambda t_0t_1+2\lambda^2-2\big)s
\\
\qquad {}+t_0^2+t_1^2+t_0t_1\big(\lambda-\lambda^{-1}\big)+2\lambda^2-2\lambda^{-2}=0.
\end{gather*}
This is a family of cubics parametrized by $(t_0,t_1,\lambda)$.
\begin{Proposition}\label{isoT} Let $\chi^*(\Lambda_0):=\{(M_0,U_1,U_2,\widehat{M}),\, u_1u_2\neq 0,\, \lambda\neq \pm 1\}/C(\mathcal{T}_0)$.
The map
$T\colon \chi(\Lambda_0)\rightarrow \mathbb{C}^6$ defined by the $6$ coordinates $(\lambda,t_0,t_1,s,x,y)$ is an isomorphism between $\chi^*(\Lambda_0)$ and the open set~$s\neq 2$, $\lambda\neq \pm 1$ in the affine variety~$F$.
\end{Proposition}

\begin{proof}
Clearly, from the Fricke lemma,
$T$ takes its values in~$F$. In order to check that this map is invertible,
we have to solve: $T(M_0,U_1,U_2,\widehat{M})=(\lambda,t_0,t_1,s,x,y)$.
This equation is equivalent to the system
\begin{subequations}\label{e(1)-e(5)}
\begin{gather}
      a_0+d_0 =t_0, \label{e(1)}\\
      \lambda a_0(1+u_1u_2)+ \lambda u_2b_0 +  \lambda^{-1}u_1c_0 +\lambda^{-1}d_0 =t_1, \label{e(2)}\\
      2+u_1u_2=s, \label{e(3)}\\
      t_0+a_0 u_1u_2+u_2b_0+u_1c_0 =x, \label{e(4)}\\
      \lambda a_0+\lambda^{-1}d_0 =y. \label{e(5)}
    \end{gather}
    \end{subequations}
If $\lambda\neq\pm 1$, we obtain from equations~\eqref{e(1)} and~\eqref{e(5)} a unique solution for $a_0$ and $d_0$
\begin{gather*}
a_0=\frac{y-\lambda^{-1}t_0}{\lambda-\lambda^{-1}},\qquad d_0=\frac{\lambda t_0-y}{\lambda-\lambda^{-1}}.
\end{gather*}
From equations~\eqref{e(2)} and~\eqref{e(4)} we obtain a unique solution for $u_1c_0$ and $u_2b_0$.
 Equation~\eqref{e(3)} gives $u_1u_2=s-2$ and we obtain a unique solution for $(\lambda,u_1u_2,u_1c_0,u_2b_0,a_0, d_0)$, which def\/ines a unique solution $(M_0,U,V,\widehat{M})$ up to the action of $C(\mathcal{T}_0)$ according to Lemma~\ref{equivalence}. Note that the solution of the system~\eqref{e(1)-e(5)} is polynomial in the variables $t_0$, $t_1$, $s$, $x$, $y$ and rational in~$\lambda$ with poles on $\lambda=\lambda^{-1}$.
\end{proof}

For $\Lambda'_0=-\Lambda_0$, we have a similar description of $\chi(\Lambda'_0)$ as an af\/f\/ine algebraic variety in the coordinates $(\lambda',t'_0,t'_1,s',x',y')$.
\begin{Proposition} The change of variables between these two charts describing the character variety~$\chi$ is given by
\begin{gather*}
      \lambda' =\lambda^{-1},\qquad       s' =s,\\
      x' =x+t_0\frac{\lambda+\lambda^{-1}}{\lambda-\lambda^{-1}}s-\frac{2}{\lambda-\lambda^{-1}}ys +
      \frac{4}{\lambda-\lambda^{-1}}y -2t_0\frac{\big(\lambda+\lambda^{-1}\big)}{\lambda-\lambda^{-1}},\\
      y' =y,
    \end{gather*}
\end{Proposition}

\begin{proof} Since from Remark~\ref{rho'} we have
\begin{gather*}
\big(M_0',U_1',U_2',\widehat{M}'\big)=P_w^{-1}\big(M_0,U_{w(1)},U_{w(2)},\widehat{M}\big)P_w=
\big(\big(M_0^{-1}\big)^t,U_2^t,U_1^t,\widehat{M}^{-1}\big),
\end{gather*} we immediately obtain
$\lambda'=\lambda^{-1}$, $s'=s$ and $y'=y$. The only non trivial computation is for $x'=\operatorname{tr}(M_0U_2U_1)$. The Fricke lemma applied on the three matrices $U_1$, $U_2$, $M_0$ gives us
\begin{gather*}
x+x'=t_0s+2b_0u_2+2c_0u_1.
\end{gather*}
From the system \eqref{e(1)-e(5)} written in the proof of Proposition~\ref{isoT}, we have $u_2b_0+u_1c_0=x-t_0-a_0(s-2)$ and $a_0=\frac{y-\lambda^{-1}t_0}{\lambda-\lambda^{-1}}$. By substitution in the above expression of~$x'$, we obtain the result.
\end{proof}

We will see in the next section that the natural parameter here is not $\Lambda_0$ but the conjugacy class~$[A_0]$ of~$A_0$ (see also   \cite[Remarks~(8.5) and~(10.6)]{Bo4}). Therefore the character variety~$\chi([A_0])$ appears as a~gluing of the two af\/f\/ine algebraic charts $\chi(\Lambda_0)$ and~$\chi(-\Lambda_0)$ which turns out to be a~scheme on our example, from the above proposition.
We can conjecture that it is the general case.

\subsection{The dynamics on the wild character variety}

This part is rather an experimental one. The way used here to encode all the information of irregular representations is not completely well justif\/ied. We need other experiments (for example in the~$P_{\rm II}$ context) and a very precise presentation of the background to present a~def\/initive version. Nevertheless, we think that a lot of tools are present here to deal with the program presented in the f\/irst section.

In the classical previous case, the conf\/iguration space $\mathcal{C}$ was the set of the positions of the 4 singularities, identif\/ied to $\mathbb{P}_1{\setminus}\{0,1,\infty\}$ through the cross ratio. In the present context we can't move the three singularities in $\mathbb{P}^1$. We can only move~$\Lambda_0$ in~$\mathcal{T}_0^*$, which is the usual conf\/iguration space considered for example in~\cite{JMU}. Ph.~Boalch in~\cite{Bo2} has introduced a coordinate independent version of this space ``through the notion of irregular curve''. Its fundamental group is again a~pure braid group.

Nevertheless, the most natural conf\/iguration space $\mathcal{C}$ is not the set of the~$\Lambda_0$'s but
 \begin{Definition}
 The conf\/iguration space $\mathcal{C}$ is the set of the initial terms~$A_0$ up to a conjugation.
 \end{Definition}

 Indeed, a gauge equivalence on $\Delta$ acts by conjugacy on the initial coef\/f\/icient~$A_0$.
In the same way, Ph.~Boalch has introduced a~``bare irregular curve'' in~\cite{Bo4}. We begin with a description of~$\mathcal{C}$ in a general framework.

Let $G$ be a reductive algebraic group, with Lie-algebra $\mathcal{G}$.
Let \begin{gather*}
\mathcal{D}=\cup\{\mathcal{T},\, \mathcal{T} \mbox{ Cartan subalgebra of } \mathcal{G}\}.
\end{gather*}
Let $\mathcal{T}_0$ be a f\/ixed Cartan subalgebra (for ${\mathfrak{sl}}_2$, $\mathcal{D}$ is the set of the \emph{diagonalisable} matrices, and~$\mathcal{T}_0$ the Cartan subalgebra of the \emph{diagonal} matrices).
Let~$\mathcal{T}^{\rm reg}$ be the subset of~$\mathcal{T}$ of matrices with distinct eigenvalues.
Recall that all the Cartan subalgebras are conjugated by some element~$g$ in~$G$. This element is not unique: if~$g$ and~$g'$ are two conjugations between~$\mathcal{T}_0$ and~$\mathcal{T}$, $g'^{-1}g$ keeps (globally) invariant~$\mathcal{T}_0$ and therefore belongs to the normalizer~$N(T_0)$ of the Cartan torus~$T_0$. We obtain a~f\/ibration
\begin{gather*}
\mathcal{D}\longrightarrow G/N(T_0),
\end{gather*}
which sends $\mathcal{T}_0$ on the identity element $I$. Note that the quotient space $G/N(T_0)$ of left classes modulo $N(T_0)$ is not a group.

\begin{Lemma} If the algebraic group $G$ is connected, and simply connected $($which is the case for ${\rm SL}_n)$, the fundamental group $\pi_1(G/N(T_0),I)$ is the Weyl group  $W:=N(T_0)/T_0$.
\end{Lemma}
\begin{proof} From the exact sequence of topological spaces
\begin{gather*}
0\rightarrow N(T_0)\rightarrow G\rightarrow G/N(T_0)\rightarrow 0
\end{gather*}
since $\pi_1(G,I)=\pi_0(G)=0$ we obtain
\begin{gather*}
0\rightarrow \pi_1(G/N(T_0),I)\rightarrow \pi_0(N(T_0))\rightarrow 0.
\end{gather*}
Now from the exact sequence of groups
\begin{gather*}
0\rightarrow T_0\rightarrow N(T_0) \rightarrow N(T_0)/T_0\rightarrow 0
\end{gather*}
since $\pi_0(T_0)=0$ we obtain $\pi_0(N(T_0))=\pi_0(N(T_0)/T_0)=W.$
\end{proof}

Remark that the f\/ibration $\mathcal{D}\longrightarrow G/N(T_0)$ has natural local trivialisations: we can lift a path~$g(s)$ in~$G/N(T_0)$ from~$A_0$ by using the conjugation~$g(s)A_0g(s)^{-1}$. The conjugation class~$[A_0]$ of~$A_0$ in~$\mathcal{D}$  cuts the f\/iber~$\mathcal{T}_0$ in a discrete set which is the orbit of~$\Lambda_0$ under the action of~$W$ by the monodromy of the above f\/ibration. Therefore, an equivalent def\/inition of the conf\/iguration space is

\begin{Definition} If the algebraic group~$G$ is connected, and simply connected,
the conf\/iguration space $\mathcal{C}$ is the set of the orbits of the Weyl group~$W$ acting on $\mathcal{T}_0^{\rm reg}$.
\end{Definition}

For $G={\rm SL}_2$, the action of the generator $w$ (the ``Weyl loop'') of $\pi_1(G/N(T_0),I)=\{\operatorname{id},w\}$ can be explicited. Consider the path in ${\rm SL}_2$ joining~$I$ to~$P_w$ by the real rotations
\begin{gather*}
w(s)= \left(\begin{matrix}
\cos s\pi/2 &-\sin s\pi/2\\ \sin s\pi/2 & \cos s\pi/2\end{matrix}\right),\qquad  s\in[0,1].
\end{gather*}
It induces a loop $w$ in $G/N(T_0)$ since $P$ belongs to $N(T_0)$. This loop can be lifted in a~path~$\delta_w$ in~$\mathcal{D}$ joining $\Lambda_0$ until $-\Lambda_0$
\begin{gather*}
\delta_w(s)=w(s)\Lambda_0w(s)^{-1}.
\end{gather*}
This path permutes the eigenvalues of~$\Lambda_0$ without moving the eigenvalues themselves.

\begin{Corollary} The whole braid group acts on~$\chi(\Lambda_0)$.
\end{Corollary}

\begin{proof}
Clearly the pure braid group $P_n=\pi_1(\mathcal{T}_0^{\rm reg},\Lambda_0)$ acts on $\chi(\Lambda_0)$, keeping invariant each element of the orbit $W\cdot \Lambda_0$.

Now let $B_n$ be the whole braid group. From the exact sequence
\begin{gather*}
0\longrightarrow P_n\longrightarrow B_n\stackrel{p}{\longrightarrow } W\longrightarrow 0,
\end{gather*}
a braid $b$ sends $\Lambda_0$ on $p(b)\cdot\Lambda_0$. Then we can lift the Weyl loop $w=p(b)$ in order to come back to~$\Lambda_0$.
\end{proof}

\looseness=-1
\textbf{The dynamics of the pure braid group.} We come back to the context of $G={\rm SL}(2,\mathbb{C})$. Let~$b$ be the (non pure) positive braid which interchanges~$\Lambda_0$ and~$-\Lambda_0$, and~$b^2$ the corresponding pure braid. The braid moves all the objects attached to~$\Lambda_0$, keeping f\/ixed~$s_\infty$. Therefore we have

\begin{Proposition}
The pure braid $b^2$ induces an automorphism $h_{b^2}$ of $\pi_1(X,S(\Lambda_0))$ which satisfies
$h_{b^2}(r_\infty)=(\gamma_{\infty,\infty})^{-1}r_\infty$, and keep invariant all the others generating morphisms.
\end{Proposition}

Now, since $\rho\circ h_{b^2}(r_\infty)=(U_1U_2\widehat{M})^{-1}$, we normalize $\rho\circ h_{b^2}$ in the equivalent representation~$\rho'$ by changing the representation of the object $s_\infty$, setting:
$X'_\infty=X_\infty\cdot (U_1U_2\widehat{M})$. Then, we set $X'_0=X_0\cdot(U_1U_2\widehat{M})$ and $X'_1=X_1\cdot(U_1U_2\widehat{M})$ in order to obtain $\rho'(\gamma_{0,\infty})=\rho'(\gamma_{1,\infty})=\rho'(\gamma_{0,1})=I$. The representation $\rho'$ is characterized by the 4-uple of matrices
\begin{gather*}
\big(M_0',U'_1,U_2',\widehat{M}'\big)=\big(\big(U_1U_2\widehat{M}\big)^{-1}M_0\big(U_1U_2\widehat{M}\big),U_1,U_2,\widehat{M}\big).
\end{gather*}

\begin{Proposition}
The automorphism $h_{b^2}$ of $\chi(\Lambda_0)$ fixes $\lambda$, $t_0$ and $t_1$, and is given in the trace coordinates on the fiber by
\begin{gather*}
      s'=s,\\
      x'=(\lambda^{-1}-\lambda+\lambda s)^2x+(\lambda^{-1}-\lambda+\lambda s)y-x\\
     \hphantom{x'=}{} -(\lambda^{-1}-\lambda+\lambda s)(\operatorname{st}_1+\lambda t_0+\lambda^{-1}t_0)+(\operatorname{st}_0+\lambda t_1+\lambda^{-1}t_1),\\
      y'=-(\lambda^{-1}-\lambda+\lambda s)x-y+(\operatorname{st}_1+\lambda t_0+\lambda^{-1}t_0).
    \end{gather*}
\end{Proposition}

\begin{proof}
We set $S:=U_1U_2$, $S':=U_1'U_2'$. Clearly $s'=\operatorname{tr}(S')=\operatorname{tr}(S)=s$. Now we have
\begin{gather*}
x'=\operatorname{tr}(M'_0S')=\operatorname{tr}\big(\widehat{M}^{-1}S^{-1}M_0S\widehat{M}S\big)
=\operatorname{tr}\big(S\widehat{M}S\widehat{M}^{-1}S^{-1}M_0\big),\\
y'=\operatorname{tr}\big(M'_0\widehat{M}'\big)=\operatorname{tr}\big(\widehat{M}^{-1}S^{-1}M_0S\widehat{M}\widehat{M}\big)
=\operatorname{tr}\big(S^{-1}M_0S\widehat{M}\big).
\end{gather*}
We can apply the extended Fricke Lemma \ref{Fricke.bis} for $M_1:=S$, $M_2:=\widehat{M}$
and $M_3:=M_0$, i.e., for $a_1=s$, $a_2=\lambda+\lambda^{-1}$, $a_3=t_0$, $a_4=t_1$, $x_{1,2}=\operatorname{tr}(S\widehat{M})=\lambda+\lambda^{-1}+\lambda(s-2)$, $x_{1,3}=x$, and $x_{2,3}=y$, which gives the expressions in the statement.
\end{proof}

\textbf{The dynamics induced by the whole braid group.} The morphism of groupoids $h_b\colon \pi_1(X,S(\Lambda_0))\rightarrow \pi_1(X,S(-\Lambda_0))$
exchanges the opposite objects related to $\Lambda_0$, f\/ixing $s_\infty$. Therefore we have
\begin{gather*}
h_b(\operatorname{st}_1)=\operatorname{st}_2,\qquad h_b(\operatorname{st}_2)=\operatorname{st}_1,\qquad h_b(r_\infty)=\big(\widehat{\beta}_{1^-}r_{1^-}\alpha_1 \big(r_{1^+}\big)^{-1}
\widehat{\beta}_{1^+}\big)^{-1}r_\infty.
\end{gather*}
Since $\rho(h_b)(r_\infty)=U_1^{-1}$ and $\rho(h_b)(\widehat{\beta}_{1^+})=\widehat{M}$, we normalize $\rho(h_b)$ in $\rho'$ by setting
\begin{gather*}
X'(s_\infty):=X(s_\infty) U_1,\qquad X'(\widehat{\tau}_2):=X(\widehat{\tau}_2)\widehat{M}^{-1},\\ X'(\widehat{\sigma}_{2^-}):=X(\widehat{\sigma}_{2^-})\widehat{M}^{-1},\qquad X'(\widehat{\sigma}_{2^+}):=X(\widehat{\sigma}_{2^+})\widehat{M}^{-1}.
\end{gather*}
The normalized representation~$\rho'$ is characterized by the 4-uple of matrices
\begin{gather*}
\big(U_1^{-1}M_0U_1,U_2,\widehat{M}U_1\widehat{M}^{-1},\widehat{M}\big).
\end{gather*}
We can check that the iteration of this action gives us
\begin{gather*}
\big(U_2^{-1}U_1^{-1}M_0U_1U_2,\widehat{M}U_1\widehat{M}^{-1},
\widehat{M}U_2\widehat{M}^{-1},\widehat{M}\big),
\end{gather*}
which is equivalent to
\begin{gather*}
\big(\big(U_1U_2\widehat{M}\big)^{-1}M_0\big(U_1U_2\widehat{M}\big),U_1,U_2,\widehat{M}\big)
\end{gather*}
and corresponds to the expression obtained above for~$h_{b^2}$. Now, if we want to obtain the expression of~$h_b$ in trace coordinates, we have to compose with the change of charts from $\chi(\Lambda_0)$ to $\chi(-\Lambda_0)$, which interchanges~$U_1$ and~$U_2$ and compose the 4-uple of matrices by~$P_w$ (see Remark~\ref{isoT}). We obtain the following action
\begin{gather*}
\big(M_0,U_1,U_2,\widehat{M}\big)\rightarrow \big(P_w\big(U_1^{-1}M_0U_1\big)P_w^{-1}, P_w\big(\widehat{M}U_1\widehat{M}^{-1}\big)P_w^{-1},P_wU_2P_w^{-1},P_w\widehat{M}P_w^{-1}\big).
\end{gather*}
In the local coordinates on~$\chi(\Lambda_0)$
\begin{gather*}
    t'_0 =t_0,\qquad
    t'_1 =t_1,\qquad
    \lambda' =\lambda^{-1},\\
    s' =\operatorname{tr}\big(P_w\big(\widehat{M}U_1\widehat{M}^{-1}U_2\big)P_w^{-1}\big)=2+\lambda^2(s-2),\\
    x' =\operatorname{tr}\big(P_w\big(U_1^{-1}M_0U_1\widehat{M}U_1\widehat{M}^{-1}U_2\big)P_w^{-1}\big)\\
\hphantom{x'}{} =   t_0+\lambda^2a_0(s-2)-\lambda^2(s-2)(u_1c_0)+(s-2)(a_0-d_0)-(s-2)(u_1c_0)+(u_2b_0),\\
        y' =\operatorname{tr}\big( P_w\big(U_1^{-1}M_0U_1\widehat{M}\big)P_w^{-1}\big) =\lambda a_0-\big(\lambda-\lambda^{-1}\big)(u_1c_0)+d_0\lambda^{-1}.
    \end{gather*}
By substituting $a_0$, $d_0$, $u_1c_0$ and $u_2b_0$ by their expressions in $\lambda$, $t_0$, $t_1$, $s$, $x$, $y$ obtained by solving the system \eqref{e(1)-e(5)} in the proof of Proposition~\ref{isoT}, we obtain explicit expressions of the dynamics which are polynomials in the variables~$t_0$, $t_1$, $s$, $x$, $y$ and rational in~$\lambda$ with poles on $\lambda-\lambda^{-1}=0$.

\section{Conclusion}

In the case of Painlev\'e~VI equation one considers a f\/ibre bundle whose base is a space of conf\/igurations and the f\/ibers character varieties. Then isomonodromy can be interpreted as an Ehresmann connection on this bundle and the fundamental group of the base (a pure braid group), induces an algebraic and symplectic dynamics on each f\/iber. Via the Riemann--Hilbert correspondence this dynamics corresponds to the dynamics of $P_{\rm VI}$.

Boalch extended partially this picture to the case of deformations of (regular-singular or irregular) $G$-connections. Introducing a f\/ibre bundle whose base is a (generalized) space of conf\/igurations and the f\/ibers (wild) character varieties, he interpreted the isomonodromy as an Ehresmann connection on this bundle. Then the fundamental group of the base (a pure braid group or a braid group) induces a symplectic dynamics on each f\/iber\footnote{One can prove that this dynamics is algebraic in the pure braids case and rational in the braids case as it is proved here for the~$P_{\rm V}$ case. We will return to such questions in a future paper.}.

In this paper we detailed the $P_{\rm V}$ case, we computed the character varieties\footnote{Using an approach based on Fricke coordinates slightly dif\/ferent of anterior computations by other authors, cf.~\cite{VdPSai}.} and the Boalch dynamics on them.
It seems that this dynamics is ``too poor'' to ref\/lect the actual ``wild dy\-na\-mics'' of~$P_{\rm V}$ (via the wild Riemann--Hilbert correspondence).
The pure braid action corresponds to the branching of the Painlev\'e functions at the f\/ixed critical points (or holonomy of the Okamoto--Painlev\'e foliation). But the interpretation of the whole braid group action remains mysterious for us.

\subsection*{Acknowledgements}
We thanks the referees for suggesting numerous improvements of a f\/irst version of this paper.
The second author thanks the \emph{Japan Society for Promotion of Science} (JSPS Fellowship for Research in Japan number S-14127) for his support for a research stay in Japan during November~2014, allowing him to present our results at the conference \emph{Recent developments in differential equations in the complex domain} (RIMS, November 17--24, 2014) and to discuss the subjects presented in this paper with some Japanese colleagues, in particular with Yousuke Ohyama and Hidetaka Sakai.

\pdfbookmark[1]{References}{ref}
\LastPageEnding

\end{document}